\documentclass[11pt]{amsart}

\usepackage{amsmath}
\usepackage{amssymb}

\newtheorem{theorem}{Theorem}[section]
\newtheorem{claim}[theorem]{Claim}

\newtheorem{corollary}[theorem]{Corollary}

\theoremstyle{definition}
\newtheorem{definition}[theorem]{Definition}

\newtheorem{question}[theorem]{Question}

\theoremstyle{remark}

\newcount\skewfactor
\def\mathunderaccent#1#2 {\let\theaccent#1\skewfactor#2
\mathpalette\putaccentunder}
\def\putaccentunder#1#2{\oalign{$#1#2$\crcr\hidewidth
\vbox to.2ex{\hbox{$#1\skew\skewfactor\theaccent{}$}\vss}\hidewidth}}

\def\smallbox#1{\leavevmode\thinspace\hbox{\vrule\vtop{\vbox
   {\hrule\kern1pt\hbox{\vphantom{\tt/}\thinspace{\tt#1}\thinspace}}
   \kern1pt\hrule}\vrule}\thinspace}

\newcommand{\cf}{{\rm cf}}

\def\qedref#1{$\qed_{\reforiginal{#1}}$}

\title{Bipartite graphs and monochromatic squares}
\author{Shimon Garti}
\address{Institute of Mathematics,
 The Hebrew University of Jerusalem,
 Jerusalem 91904, Israel}
\email{shimon.garty@mail.huji.ac.il}

\begin{document}
\let\labeloriginal\label
\let\reforiginal\ref

\begin{abstract}
Let $\kappa$ be a successor cardinal.
We prove that consistently every bipartite graph of size $\kappa^{+}\times\kappa^{+}$ contains either an independent set or a clique of size $\tau\times\tau$ for every ordinal $\tau<\kappa^{+}$. We prove a similar theorem for $\ell$-partite graphs.
\end{abstract}

\maketitle

\newpage

\section{Introduction}

Let $G=(U\times V,E)$ be a bipartite graph.
An independent set in $G$ is a product $A\times B$ such that $A\subseteq U, B\subseteq V$ and $(\alpha,\beta)\notin E$ whenever $\alpha\in A, \beta\in B$.
A complete subgraph (or a clique) is a product $A\times B$ such that $A\subseteq U, B\subseteq V$ and $(\alpha,\beta)\in E$ for every $\alpha\in A, \beta\in B$.
If one wishes to eliminate large independent sets then many edges must be added everywhere, so the graph is likely to contain large complete subgraphs.
The following \emph{mot d'ordre} obeys this basic intuition:

\begin{center}
Every bipartite graph contains either \newline a large independent set or a large clique.
\end{center}

The natural question, therefore, is how large?
The following claim gives a simple limitation.
We indicate that our bipartite graphs are balanced in the sense that $|U|=|V|$ where $G=(U\times V,E)$.

\begin{claim}
\label{clmlim} Let $\kappa$ be an infinite cardinal. \newline
There exists a $\kappa\times\kappa$ bipartite graph with no independent set of size $\kappa\times\kappa$ and no complete subgraph of size $\kappa\times\kappa$.
\end{claim}

\par\noindent\emph{Proof}. \newline
Let $G=(U\times V,E)$ where $U=V=\kappa$ and $(\alpha,\beta)\in E$ iff $\alpha\leq\beta$.
Suppose that $A,B\subseteq\kappa$ and $|A|=|B|=\kappa$.
Choose any $\alpha\in A$ and choose $\beta\in B$ such that $\beta\geq\alpha$ (possible, since $|B|=\kappa$).
By definition $(\alpha,\beta)\in E$, and hence $A\times B$ is not independent.
Now choose $\gamma\in A$ such that $\gamma>\beta$ (possible, since $|A|=\kappa$).
Again, by the definition of the graph we see that $(\gamma,\beta)\notin E$ so $A\times B$ is not a complete subgraph.

\hfill \qedref{clmlim}

Actually, the proof gives more than stated.
It shows that one can always produce a graph with no independent set of size $1\times\kappa$ and no clique of size $\kappa\times 1$.
In this light, the main result of this paper is optimal (when balanced products are considered):

\begin{theorem}
\label{thmmt} Let $\kappa$ be a successor cardinal. \newline
It is consistent that for every bipartite graph $G$ of size $\kappa^{+}\times\kappa^{+}$ and every ordinal $\tau\in\kappa^{+}$, the graph $G$ contains either an independent set or a complete subgraph of order type $\tau\times\tau$.
\end{theorem}

\hfill \qedref{thmmt}

Let us mention some concepts that will be used throughout the paper, and fix our notation. The ambient combinatorial notion behind the results to follow is called the polarzied partition relation. We write $\binom{\alpha}{\beta} \rightarrow \binom{\gamma}{\delta}_\theta$ iff for every coloring $c:\alpha\times\beta\rightarrow\theta$ one can find $A\subseteq\alpha, B\subseteq\beta$ so that ${\rm otp}(A)=\gamma, {\rm otp}(B)=\delta$ and $c\upharpoonright(A\times B)$ is monochromatic.

If $\kappa=\cf(\kappa)<\lambda$ then $S^\lambda_\kappa$ is the set $\{\delta\in\lambda:\cf(\delta)=\kappa\}$. If $\cf(\lambda)>\omega$ then this set is stationary in $\lambda$.
An elementary embedding $\jmath:V\prec M$ is a non-trivial embedding of the universe of set theory into a transitive model $M$ of set theory.
A cardinal $\kappa$ is the critical point of $\jmath$ iff $\kappa$ is the first ordinal moved by $\jmath$. This fact is denoted by $\kappa={\rm crit}(\jmath)$.
We call $\kappa$ a huge cardinal iff there exists an elementary embedding $\jmath:V\prec M$ so that $\kappa={\rm crit}(\jmath)$ and ${}^{\jmath(\kappa)}M\subseteq M$.
For a general background regarding large cardinals we suggest \cite{MR1994835}.

Let $\mathcal{I}$ be an ideal over $\kappa$.
We denote the collection of all subsets of $\kappa$ which are not in $\mathcal{I}$ by $\mathcal{I}^+$. We shall call these sets $\mathcal{I}$-positive sets.
An ideal $\mathcal{I}$ is $(\lambda,\mu,\theta)$-saturated iff for every collection $\{A_\alpha:\alpha\in\lambda\}\subseteq\mathcal{I}^+$ one can find a subcollection $\{A_{\alpha_\beta}:\beta\in\mu\}$ such that $\bigcap\{A_{\alpha_\beta}:\beta\in\mathcal{C}\}$ is $\mathcal{I}$-positive whenever $\mathcal{C}\in[\mu]^\theta$.

A directed graph $G=(V,E)$ is called a directed bipartite graph iff there exists a partition of $V$ into two disjoint sets $U,W$ such that each of them is edge-free.
We use the ordered pair notation $\langle\alpha,\beta\rangle\in E$ to indicate that $E$ contains an edge from $\alpha$ to $\beta$.
A multigraph is a triple $G=(V,E,\psi)$ where $V$ is a set of vertices, $E$ is a set of edges and $\psi:E\rightarrow\{(x,y):x,y\in V\}$ is a function which spells out the endpoints of each edge in $E$.
The idea is that we allow parallel edges between the same two vertices, so we specify with $\psi$ the pair of vertices for which some $e\in E$ is applied to.
We adopt the convention that multigraphs are loopless.
The concept of a bipartite multigraph is defined similarly.
As a general reference regarding concepts in graph theory we suggest \cite{MR2744811}.

The rest of the paper contains three additional sections.
In the first section we focus on square in square with $3$ colors.
The main result is based on a theorem from \cite{MR1833480}, and it gives a monochromatic square of order type $\kappa+1\times\kappa+1$.
The derived conclusions apply to bipartite graphs and directed bipartite graphs.
In the second section we prove the consistency of stronger combinatorial relations with infinitely many colors.
This is done for double successors of regular cardinals.
A large number of colors is useful when dealing with multigraphs.
We generalize the results about bipartite graphs to $\ell$-partite graphs.
In the third section we point to some open problems.

\newpage

\section{Bipartite graphs}

The combinatorial tool that we shall apply to bipartite graphs is called square in square.
The idea is that we find a monochromatic square of size $\tau\times\tau$ for every coloring defined on a larger square of size $\kappa\times\kappa$.
This concept is a special case of the polarized relation.

\begin{definition}
\label{defsis} Square in square. \newline
Assume that $\theta,\kappa$ are cardinals, $\theta<\kappa$ and $\tau\in\kappa$. \newline
We say that ${\rm sis}(\kappa,\tau,\theta)$ holds iff for every coloring $c:\kappa\times\kappa\rightarrow\theta$ there are $A,B\subseteq\kappa, {\rm otp}(A)={\rm otp}(B)=\tau$ such that $c\upharpoonright(A\times B)$ is constant.
\end{definition}

In the terminology of polarized relations, ${\rm sis}(\kappa,\tau,\theta)$ is the positive relation $\binom{\kappa}{\kappa}\rightarrow \binom{\tau}{\tau}^{1,1}_\theta$.
The following basic claim demonstrates the way in which theorems about bipartite graphs can be deduced from instances of square in square with two colors.

\begin{claim}
\label{clmsis} Assume ${\rm sis}(\kappa,\tau,2)$. \newline
Then every bipartite graph of size $\kappa\times\kappa$ contains either an independent set or a complete subgraph of order type $\tau\times\tau$.
\end{claim}

\par\noindent\emph{Proof}. \newline
Let $G=(U\times V,E)$ be a bipartite graph of size $\kappa\times\kappa$.
Define $c:U\times V\rightarrow 2$ by $c(\alpha,\beta)=0$ iff $(\alpha,\beta)\notin E$.
By the assumption ${\rm sis}(\kappa,\tau,2)$ there are $A\subseteq U, B\subseteq V, {\rm otp}(A)={\rm otp}(B)=\tau$ such that $A\times B$ is $c$-monochromatic.
Now if $c''(A\times B)=\{0\}$ then $(\alpha,\beta)\notin E$ whenever $\alpha\in A,\beta\in B$ so $A\times B$ is independent.
If $c''(A\times B)=\{1\}$ then $(\alpha,\beta)\in E$ whenever $\alpha\in A,\beta\in B$ and hence $A\times B$ is a clique.

\hfill \qedref{clmsis}

The main result of this section follows from a theorem of Baumgartner and Hajnal which applies to colorings with three colors.
In the context of cliques and independent sets, the effect of three colors is best demonstrated in directed bipartite graphs.
But we have to be accurate about the definition of a complete subgraph.
Suppose that $G=(U\times V,E)$ is a directed bipartite graph and $A\subseteq U, B\subseteq V$. Let $H$ be the induced subgraph.

If we define $H$ to be complete iff both $\langle\alpha,\beta\rangle$ and $\langle\beta,\alpha\rangle$ belong to $E$ for every $\alpha\in A, \beta\in B$ then it is easy to construct a directed bipartite graph with no independent subgraph and no complete subgraph which are not empty.
Indeed, let $E=\{\langle\alpha,\beta\rangle:\alpha\in U,\beta\in V\}$.
In order to enable non-empty statements we define $H$ to be a complete subgraph iff for every $\alpha\in A,\beta\in B$ either $\langle\alpha,\beta\rangle\in E$ or $\langle\beta,\alpha\rangle\in E$.

\begin{claim}
\label{clmdibipartite} Assume ${\rm sis}(\kappa,\tau,3)$. \newline
Then every $\kappa\times\kappa$ directed bipartite graph contains either an independent set or a complete subgraph of order type $\tau\times\tau$.
\end{claim}

\par\noindent\emph{Proof}. \newline 
Let $G=(U\times V,E)$ be a directed bipartite graph.
Define a coloring $c:U\times V\rightarrow 3$ as follows:
$$
c(\alpha,\beta)=
\begin{cases}
0 & \text{if } \langle\alpha,\beta\rangle\in E\\
1 & \text{if } \langle\alpha,\beta\rangle\notin E \wedge \langle\beta,\alpha\rangle\in E\\
2 & \text{if } \langle\alpha,\beta\rangle\notin E \wedge \langle\beta,\alpha\rangle\notin E
\end{cases}
$$
By ${\rm sis}(\kappa,\tau,3)$ there are $A\subseteq U, B\subseteq V$ and $i<3$ such that ${\rm otp}(A)={\rm otp}(B)=\tau$ and $c''(A\times B)=\{i\}$.
The induced directed subgraph $H$ will be a complete $\tau\times\tau$ subgraph when $i\in\{0,1\}$ and a $\tau\times\tau$ independent subset when $i=2$, so we are done.

\hfill \qedref{clmdibipartite}

We quote now Theorem 2.1 from \cite{MR1833480}. The statement there gives a monochromatic $\kappa\times\kappa$ square, but the proof shows that actually one has a square of size $\kappa+1\times\kappa+1$.

\begin{theorem}[Baumgartner-Hajnal]
\label{thmmt1} Assume that $\kappa=\kappa^{<\kappa}$. \newline
Then ${\rm sis}(\kappa^+,\kappa+1,3)$.
\end{theorem}

\hfill \qedref{thmmt1}

We can derive now the following conclusion:

\begin{corollary}
\label{cor1} Let $G$ be any $\kappa^+\times\kappa^+$ bipartite graph, where $\kappa$ is a regular cardinal.
\begin{enumerate}
\item [$(\aleph)$] If $\kappa$ is strongly inaccessible then there is a $\kappa+1\times\kappa+1$ subgraph of $G$ which is either complete or independent.
\item [$(\beth)$] If $2^\kappa=\kappa^+$ then there is a $\kappa+1\times\kappa+1$ subgraph of $G$ which is either complete or independent.
\end{enumerate}
\end{corollary}

We do not know whether ${\rm sis}(\kappa^+,\tau,3)$ can be proved for every $\kappa=\kappa^{<\kappa}$ and every $\tau\in\kappa^+$.
In many cases this can be forced, as will be shown in the next section.
We shall discuss this issue in the last section.

\newpage 

\section{Many colors and $\ell$-partite graphs}

In this section we show how to force instances of square in square with many colors. As mentioned above, there are many limitations on this possibility as one can force a negative result (in many cases) even when the number of colors is just $4$, see \cite{MR1833480}.
Our ability to force the positive side of the coin concentrates on double successors of cardinals.
This is done here for double successors of regular cardinals.

We shall assume that there is a huge cardinal above $\kappa$. 
We modify the proof of Theorem 3.2 from \cite{MR3610266}, in order to get square in square. It is not clear whether the method of \cite{MR3610266} can be applied to successors of singular cardinals, see the next section.

\begin{theorem}
\label{mt1} Double successors of regular cardinals. \newline
Assume $\theta=\cf(\theta), \kappa=\theta^+$ and there exists a huge cardinal above $\theta$. \newline 
Then one can force ${\rm sis}(\kappa^+,<\kappa^+,\theta)$.
\end{theorem}

\par\noindent\emph{Proof}. \newline
Let $\mathcal{I}$ be a $\theta^+$-complete and $(\theta^{++},\theta^{++},\theta)$-saturated ideal over $\theta^+$. This can be forced by \cite{MR673792}. Notice that $2^\theta=\theta^+$ in this construction. Fix an ordinal $\tau\in(\theta^+,\theta^{++})$, and a bijection $h:\theta^+\rightarrow\tau$.

We define an ideal $\mathcal{J}$ over $\tau$ as follows. Enumerate the elements of $\mathcal{I}$ by $\{I_\varepsilon:\varepsilon<\delta_*\}$ and let $J_\varepsilon = h''I_\varepsilon$ for every $\varepsilon<\delta_*$. Set $\mathcal{J} = \{J_\varepsilon:\varepsilon<\delta_*\}$ and notice that $\mathcal{J}$ is an ideal over $\tau$. Moreover, $\mathcal{J}$ is $\theta^+$-complete and $(\theta^{++},\theta^{++},\theta)$-saturated. This follows from the fact that $h$ is a bijection.

Let $c:\theta^{++}\times\theta^{++}\rightarrow\theta$ be any coloring. For every $\alpha<\theta^{++},\iota<\theta$ set $x^\iota_\alpha = \{\beta\in\tau: c(\alpha,\beta)=\iota\}$. Notice that we focus here on $c\upharpoonright(\theta^{++}\times\tau)$, so actually we prove that $\binom{\theta^{++}}{\tau}\rightarrow \binom{\tau}{\tau}^{1,1}_\theta$ for every $\tau<\theta^{++}$. By the $\theta^+$-completeness of $\mathcal{J}$, for each $\alpha<\theta^{++}$ there exists $\iota(\alpha)<\theta$ so that $x_\alpha^{\iota(\alpha)}\in\mathcal{J}^+$.

Since $\theta^{++}=\cf(\theta^{++})>\theta$, there is some $A\subseteq \theta^{++}, |A|=\theta^{++}$ and $\iota<\theta$ such that $\alpha\in A \Rightarrow \iota(\alpha)=\iota$. As all we need is a monochromatic product whose large component is of size $\theta^{++}$, we can restrict the coloring to $A\times\theta^{++}$ or assume without loss of generality that $A=\theta^{++}$. So let $x=\{x^\iota_\alpha:\alpha<\theta^{++}\}\subseteq \mathcal{J}^+$, and we omit $\iota$ and write $x_\alpha$ from now on.

Choose a large enough regular $\chi$ (e.g., $\chi = (\beth_3(\theta^+))^+$). Choose an elementary chain $(M_\eta:\eta\leq\tau)$ of submodels of $\mathcal{H}(\chi)$ such that for every $\eta\leq\tau$ the following requirements are met:
\begin{enumerate}
\item [$(\aleph)$] $|M_\eta|=\theta^+$ and $\tau\subseteq M_\eta$.
\item [$(\beth)$] $\mathcal{J},c,\tau,h,x\in M_\eta$.
\item [$(\gimel)$] ${}^{\leq\theta}M_\eta\subseteq M_\eta$.
\item [$(\daleth)$] $\zeta<\eta\leq\tau\Rightarrow M_\zeta\in M_\eta$.
\end{enumerate}
We comment that here we require $\tau\subseteq M_\eta$ (rather than just $\theta^+\subseteq M_\eta$ as in Theorem 3.2 of \cite{MR3610266}), so we will be able to produce two sets of order type $\tau$ in the monochromatic product.
We denote $\sup(M_\eta\cap\theta^{++})$ by $\sigma_\eta$, for every $\eta\leq\tau$. Without loss of generality, $\bigcap\{x_\alpha:\alpha\in C\}\in \mathcal{J}^+$ for any $C\in[\theta^{++}]^\theta$. Let $\delta$ be $\sigma_\tau$.
We try to define, by induction on $\gamma<\theta^+$, two sequences of ordinals. The sequences are of the form $\langle\alpha_{h(\gamma)} :\gamma<\theta^+\rangle$ and $\langle\beta_{h(\gamma)} :\gamma<\theta^+\rangle$, and all the elements are ordinals of $\tau$.

For $\gamma=0$ we define, first, $\beta_{h(0)}=\min(x_\delta)$. We choose now an ordinal $\zeta\in M_{h(0)+1}\setminus M_{h(0)}$ such that $\beta_{h(0)}\in x_\zeta$. Such $\zeta$ exists since $\delta>\sigma_{h(0)+1}, \beta_{h(0)}\in x_\delta$ and by elementarity. Let $\alpha_{h(0)}$ be $\zeta$.

Assume now that $\gamma>0$ and the construction has been completed up to this stage. Let $\beta_{h(\gamma)} = \min(\bigcap\limits_{\xi<\gamma}x_{\alpha_{h(\xi)}}\cap x_\delta \setminus \{\beta_{h(\xi)}:\xi<\gamma\})$. The saturation of the ideal assures us that $\beta_{h(\gamma)}$ is well-defined. Let $B = \{\beta_{h(\xi)}:\xi<\gamma\}$. Notice that $B\in M_\eta$ for every $\eta\leq\tau$, by requrirement $(\gimel)$ above. Moreover, $B\subseteq x_\delta$ since we added $x_\delta$ to the intersections at every step in the choice of the $\beta$-s. It follows that, for each $\eta, M_\eta\models$ there are unboundedly many $\alpha<\theta^{++}$ such that $B\subseteq x_\alpha$. Hence we can choose $\alpha_{h(\gamma)}\in M_{h(\gamma)+1}\setminus M_{h(\gamma)}$ with the property of $B\subseteq x_{\alpha_{h(\gamma)}}$.

Having accomplished the inductive process, let $A = \{\alpha_{h(\gamma)}:\gamma<\theta^+\}, B = \{\beta_{h(\gamma)}: \gamma<\theta^+\}$. Notice that ${\rm otp}(A,<) = {\rm otp}(B,<) = \tau$. Likewise, $c''(A\times B) = \{\iota\}$, so we are done.

\hfill \qedref{mt1}

We derive from the above theorem a conclusion about bipartite multigraphs.
Suppose that $G=(U\times V,E,\psi)$ is a bipartite multigraph.
Fix an enumeration $\{e_\alpha:\alpha\in\delta\}$ of the elements of $E$.
For every $\alpha\in U,\beta\in V$ let $\eta(\alpha,\beta)$ be ${\rm otp}\{e\in E:\psi(e)=(\alpha,\beta)\}$ ordered according to the above fixed enumeration.

\begin{theorem}
\label{thmmulti} Let $\kappa$ be a regular cardinal. \newline 
It is consistent (by assuming that there is a huge cardinal above $\kappa$) that for every $\kappa^{++}\times\kappa^{++}$ bipartite multigraph $G=(U\times V,E,\psi)$ which satisfies $\sup\{\eta(\alpha,\beta):\alpha\in U, \beta\in V\}\leq\kappa$ and any ordinal $\tau\in\kappa^{++}$ one can find $A\subseteq U, B\subseteq V, {\rm otp}(A)={\rm otp}(B)=\tau$ such that either $A\times B$ is independent or $A\times B$ is a clique in which $\eta(\alpha,\beta)$ assumes a constant value for every $\alpha\in A,\beta\in B$.
\end{theorem}

\par\noindent\emph{Proof}. \newline 
We may assume that $U=V=\kappa^{++}$.
Force ${\rm sis}(\kappa^{++},<\kappa^{++},\kappa)$.
Define $c:\kappa^{++}\times\kappa^{++}\rightarrow\kappa$ by $c(\alpha,\beta)=0$ when $(\alpha,\beta)\notin E$ and $c(\alpha,\beta)=\eta(\alpha,\beta)$ otherwise.
Now use ${\rm sis}(\kappa^{++},<\kappa^{++},\kappa)$ in order to derive the desired conclusion.

\hfill \qedref{thmmulti}

We conclude with a theorem about monochromatic $\ell$-partite graphs.
If $G=(U_1\times\cdots\times U_k,E)$ is a $k$-partite graph then one can find a bipartite subgraph which forms either a large clique or a large independent subset by the results of the previous section.
The interesting question is whether such a large subgraph can be $\ell$-bipartite for $\ell>2$.

By forcing intervals of cardinals which satisfy the continuum hypothesis one can try to apply instances of ${\rm sis}(\kappa^{+n},<\kappa^{+(n-1)},3)$ finitely many times with the goal of producing an $\ell$-partite subgraph of size $\kappa\times\kappa$ which is either a clique or independent.
At every step the size of the subgraph is reduced by one cardinality.
It seems that some pairs coloring (with two colors) is required in order to get the same option (a clique or an independent set) for many pairs simultaneously.
We shall use the classical theorem of Ramsey.
Denote by $R(\ell)$ the minimal $n$ such that any coloring of $[n]^2$ with two colors has a monochromatic subset of size $\ell$.

\begin{theorem}
\label{thmlpartite} Let $\kappa$ be an infinite cardinal, $\ell\in\omega$.
\newline 
Assume that $2^\lambda=\lambda^+$ for every $\lambda<\kappa^{+\omega}$. \newline 
Let $n=R(\ell), m=n-1$. \newline 
Then every $n$-partite $\kappa^{+m}\times\kappa^{+m}$ graph contains an $\ell$-partite $\kappa\times\kappa$ subgraph which forms either a clique or an independent set.
\end{theorem}

\par\noindent\emph{Proof}. \newline 
Let $G=(U_0\times\cdots\times U_{n-1},E)$ be any $n$-partite graph such that $|U_i|=\kappa^{+m}$ for every $i<n$.
For each $i<n$ let $A^i_0=U_i$.
We use double induction in order to prove the above statement.
By induction on $i<n-1$ we define $A^i_j$ for every $j\in(i,n)$ as follows.

We choose $A^i_j\subseteq A^i_{j-1}$ such that $|A^i_j|^+=|A^i_{j-1}|$ and for some $B\subseteq A^j_i, |B|=|A^i_j|$ we have either that $A^i_j\times B$ is a clique or an independent set.
Then we define $A^j_{i+1}=B$ and proceed with $j$.
After $m$ steps we create for each $i<n$ an $\subseteq$-decreasing sequence of sets $\langle A^i_k:k<m\rangle$ such that $|A^i_m|=\kappa$.
Likewise, for every $0\leq i_0<i_1\leq m$ either $A^{i_0}_m\times A^{i_1}_m$ is a $\kappa\times\kappa$ clique or a $\kappa\times\kappa$ independent subset of $G$.

Define a coloring $d:[n]^2\rightarrow\{0,1\}$ by letting $d(i_0,i_1)=0$ iff $A^{i_0}_m\times A^{i_1}_m$ is a clique, and $d(i_0,i_1)=1$ iff $A^{i_0}_m\times A^{i_1}_m$ is independent.
By Ramsey's theorem there exists a set $y\subseteq n,|y|=\ell$ such that $y$ is $d$-monochromatic.
The induced subgraph over $(A^i_m)_{i\in y}$ is as desired.

\hfill \qedref{thmlpartite}

The values of $m,n$ in the above theorem depend on $R(\ell)$.
It is easy to see that if $n<R(\ell)$ then there is a coloring of the $n$-partite graph of size $\kappa^{n-1}\times\kappa^{n-1}$ with no $\ell$-partite clique and no $\ell$-partite independent set, so Ramsey number is optimal here.
One may wonder if these results can be generalized to $\theta$-partite subgraphs of $\lambda$-partite graphs where $\theta$ and $\lambda$ are infinite.

\newpage 

\section{Open problems}

The first question that we suggest is whether the theorem of Baumgartner and Hajnal can be improved.
Remark that in many cases an improvement can be forced, but one may wonder about the ZFC situation.

\begin{question}
\label{qbaha} Assume that $\kappa=\kappa^{<\kappa}$. \newline 
Is the relation ${\rm sis}(\kappa^+,<\kappa^+,3)$ provable in ZFC?
\end{question}

'e's never 'appy unless 'e's miserable (\cite[p. 398]{herriot}).
Upon replacing `miserable' by `measurable' we indicate that if $\kappa$ is measurable then the answer to the above question is positive, using any normal ultrafilter over $\kappa$ and the main result of the previous section.

The above question is connected with a fascinating problem related to triplets.
The majority of results in the literature apply to pairs.
There is some reference to colorings with larger domains, notwithstanding.
For the next problem let us generalize the notation of square in square.

\begin{definition}
\label{defcic} Cube in cube. \newline 
Assume that $\theta,\kappa$ are cardinals, $\theta<\kappa$ and $\tau\in\kappa$. \newline 
We say that ${\rm cic}(\kappa,\tau,\theta)$ holds iff for every coloring $c:\kappa\times\kappa\times\kappa\rightarrow\theta$ there are $H,I,J\subseteq\kappa$ of order type $\tau$ such that $c\upharpoonright(H\times I\times J)$ is constant.
\end{definition}

Using this terminology, Question 28 from \cite{MR0280381} is whether ${\rm cic}(\aleph_1,\aleph_0,2)$ holds.
The problem is generalized in \cite[p. 108]{MR3075383}, and can be phrased as follows:

\begin{question}
\label{qcube} Let $\kappa$ be an infinite cardinal. 
\begin{enumerate}
\item [$(\aleph)$] Is it provable, under any assumption on $\kappa$, that ${\rm cic}(\kappa^+,\kappa,2)$ holds?
\item [$(\beth)$] Is it consistent for some $\kappa$ that ${\rm cic}(\kappa^+,\kappa,2)$ holds?
\end{enumerate}
\end{question}

It seems that no progress has been made regarding these questions, neither by proving positive statements nor by giving counterexamples.
We indicate that the methods of \cite{MR1833480} might be useful here.
We also conjecture that if $\kappa>\cf(\kappa)$ then $\neg{\rm cic}(\kappa^+,\kappa,2)$.

This brings us to the issue of singular cardinals.
The main theorem of the current paper deals with double successors of \emph{regular cardinals}. What can be sais about double successors of singular cardinals?
Basically, there are two approaches here.
The first one is to try forcing a sufficiently saturated ideal over a successor of a singular cardinal.
The second approach is to force a saturated ideal over $\kappa^+$ where $\kappa$ is measurable and then to singularize $\kappa$ with the hope that ${\rm sis}(\kappa^{++},<\kappa^{++},\theta)$ will remain in the generic extension even if the saturation fades away.

The first part of this plan is possible if there exists a supercompact cardinal $\kappa$ and a huge cardinal above $\kappa$ (remark that this assumption implies that there are many huge cardinals).
If one forces $\kappa$ to be Laver-indestructible then a very saturated ideal over $\kappa^+$ is obtainable while preserving the supercompactness of $\kappa$, since Laver's forcing to create such an ideal is $\kappa$-directed-closed.
For the second part, suppose that one adds a Prikry sequence to $\kappa$.
It is not clear whether the saturation is preserved, but maybe ${\rm sis}(\kappa^{++},<\kappa^{++},\theta)$ is easier to keep.

For this end, one has to force $2^\kappa=\kappa^+$ (this is possible to be added to the forcing which makes $\kappa$ indestructible).
Now if $\lambda=\cf(\lambda)>\kappa^+$ then any set of size $\lambda$ in the generic extension contains a set of size $\lambda$ from the ground model.
Thus any new coloring $c:\kappa^{++}\times\kappa^{++}\rightarrow\theta$ contains a coloring from the ground model of size $\kappa^{++}$.
This old coloring has a monochromatic product which will serve for $c$ as well.
The problem here is to make sure that the old coloring has size $\kappa^{++}$ in both coordinates.
This leads to the following:

\begin{question}
\label{qsing} Is it consistent that ${\rm sis}(\kappa^{++},<\kappa^{++},\theta)$ holds for some $\theta\geq 2$ and $\kappa>\cf(\kappa)$?
\end{question}

\newpage

\bibliographystyle{amsplain}
\bibliography{arlist}

\providecommand{\bysame}{\leavevmode\hbox to3em{\hrulefill}\thinspace}
\providecommand{\MR}{\relax\ifhmode\unskip\space\fi MR }
\providecommand{\MRhref}[2]{%
  \href{http://www.ams.org/mathscinet-getitem?mr=#1}{#2}
}
\providecommand{\href}[2]{#2}
\begin{thebibliography}{1}

\bibitem{MR1833480}
James~E. Baumgartner and Andras Hajnal, \emph{Polarized partition relations},
  J. Symbolic Logic \textbf{66} (2001), no.~2, 811--821. \MR{1833480}

\bibitem{MR2744811}
Reinhard Diestel, \emph{Graph theory}, fourth ed., Graduate Texts in
  Mathematics, vol. 173, Springer, Heidelberg, 2010. \MR{2744811}

\bibitem{MR0280381}
P.~Erd\H{o}s and A.~Hajnal, \emph{Unsolved problems in set theory}, Axiomatic
  {S}et {T}heory ({P}roc. {S}ympos. {P}ure {M}ath., {V}ol. {XIII}, {P}art {I},
  {U}niv. {C}alifornia, {L}os {A}ngeles, {C}alif., 1967), Amer. Math. Soc.,
  Providence, R.I., 1971, pp.~17--48. \MR{0280381}

\bibitem{MR3610266}
Shimon Garti, \emph{Amenable colorings}, Eur. J. Math. \textbf{3} (2017),
  no.~1, 77--86. \MR{3610266}

\bibitem{herriot}
James Herriot, \emph{All things bright and beautiful}, St. Martin's Paperbacks,
  1974.

\bibitem{MR1994835}
Akihiro Kanamori, \emph{The higher infinite}, second ed., Springer Monographs
  in Mathematics, Springer-Verlag, Berlin, 2003, Large cardinals in set theory
  from their beginnings. \MR{1994835 (2004f:03092)}

\bibitem{MR673792}
Richard Laver, \emph{An {$(\aleph _{2},\,\aleph _{2},\,\aleph _{0})$}-saturated
  ideal on {$\omega _{1}$}}, Logic {C}olloquium '80 ({P}rague, 1980), Stud.
  Logic Foundations Math., vol. 108, North-Holland, Amsterdam-New York, 1982,
  pp.~173--180. \MR{673792}

\bibitem{MR3075383}
Neil~H. Williams, \emph{Combinatorial set theory}, Studies in Logic and the
  Foundations of Mathematics, vol.~91, North-Holland Publishing Co., Amsterdam,
  1977. \MR{3075383}

\end{thebibliography}

\end{document}